\documentclass[11pt]{article}
\usepackage[centering,includeheadfoot,margin=2 cm]{geometry}
\usepackage{graphics,enumitem,epsfig}
\usepackage{amsfonts,amsmath,amssymb,euscript,color}
\usepackage{pst-all}

\begin{document}

\def\ALERT#1{{\large\bf $\clubsuit$#1$\clubsuit$}}

\numberwithin{equation}{section}
\newtheorem{defin}{Definition}
\newtheorem{theorem}{Theorem}
\newtheorem{proposition}{Proposition}
\newtheorem{notice}{Notice}
\newtheorem{hypothesis}{Hypothesis}
\newtheorem{lemma}{Lemma}
\newtheorem{cor}{Corollary}
\newtheorem{example}{Example}
\newtheorem{remark}{Remark}
\newtheorem{conj}{Conjecture}
\def\begproof{\noindent{\bf Proof: }}
\def\endproof{\par\rightline{\vrule height5pt width5pt depth0pt}\medskip}
\def\div{\nabla\cdot}
\def\rot{\nabla\times}
\def\sign{{\rm sign}}
\def\arsinh{{\rm arsinh}}
\def\arcosh{{\rm arcosh}}
\def\diag{{\rm diag}}
\def\const{{\rm const}}
\def\d{\,\mathrm{d}}
\def\eps{\varepsilon}
\def\phi{\varphi}
\def\theta{\vartheta}
\def\N{\mathbb{N}}
\def\R{\mathbb{R}}
\def\C{\hbox{\rlap{\kern.24em\raise.1ex\hbox
      {\vrule height1.3ex width.9pt}}C}}
\def\P{\hbox{\rlap{I}\kern.16em P}}
\def\Q{\hbox{\rlap{\kern.24em\raise.1ex\hbox
      {\vrule height1.3ex width.9pt}}Q}}
\def\M{\hbox{\rlap{I}\kern.16em\rlap{I}M}}
\def\Z{\hbox{\rlap{Z}\kern.20em Z}}
\def\({\begin{eqnarray}}
\def\){\end{eqnarray}}
\def\[{\begin{eqnarray*}}
\def\]{\end{eqnarray*}}
\def\part#1#2{\frac{\partial #1}{\partial #2}}
\def\partk#1#2#3{\frac{\partial^#3 #1}{\partial #2^#3}} 
\def\mat#1{{D #1\over Dt}}
\def\dx{\nabla_x}
\def\dv{\nabla_v}
\def\grad{\nabla}
\def\Norm#1{\left\| #1 \right\|}
\def\pmb#1{\setbox0=\hbox{$#1$}
  \kern-.025em\copy0\kern-\wd0
  \kern-.05em\copy0\kern-\wd0
  \kern-.025em\raise.0433em\box0 }
\def\bar{\overline}
\def\lbar{\underline}
\def\fref#1{(\ref{#1})}
\def\half{\frac{1}{2}}
\def\oo#1{\frac{1}{#1}}

\def\tot#1#2{\frac{\d #1}{\d #2}} 
\def\laplace{\Delta}
\def\d{\,\mathrm{d}}
\def\N{\mathbb{N}}
\def\R{\mathbb{R}}
\def\supp{\mbox{supp }}
\def\eps{\varepsilon}
\def\phi{\varphi}

\def\J{\mathcal{J}}
\def\I{\mathcal{I}}
\def\K{\mathcal{K}}

\def\Ikn{\I_{k_n}}
\def\sIkn{(\Ikn)_{n\in\N}}

\def\O{\mathcal{O}}
\def\calG{\mathcal{G}}
\def\S{\mathcal{S}}
\def\calW{\mathcal{W}}

\def\td{\alpha}

\def\comment#1{\textbf{#1}}


\centerline{{\huge Flocking dynamics and mean-field limit}}
\centerline{{\huge in the Cucker-Smale-type model}}
\centerline{{\huge with topological interactions}}


\vskip 5mm

\centerline{
{\large Jan Haskovec}\footnote{Mathematical and Computer Sciences and Engineering Division,
King Abdullah University of Science and Technology,
Thuwal 23955-6900, Kingdom of Saudi Arabia; 
{\it jan.haskovec@kaust.edu.sa}}
}
\vskip 6mm


\noindent{\bf Abstract.}
We introduce a Cucker-Smale-type model for flocking, where the strength of interaction between two agents
depends on their relative separation (called ``topological distance'' in previous works), 
which is the number of intermediate individuals separating them. This makes the model scale-free
and is motivated by recent extensive observations of starling flocks, suggesting that
interaction ruling animal collective behavior depends on topological rather than metric distance.
We study the conditions leading to asymptotic flocking in the topological model,
defined as the convergence of the agents' velocities to a common vector.
The shift from metric to topological interactions requires development of new analytical methods,
taking into account the graph-theoretical nature of the problem. Moreover, we provide a rigorous
derivation of the mean-field limit of large populations, recovering kinetic and hydrodynamic descriptions.
In particular, we introduce the novel concept of relative separation in continuum descriptions,
which is applicable to a broad variety of models of collective behavior.
As an example, we shortly discuss a topological modification of the attraction-repulsion model
and illustrate with numerical simulations that the modified model produces interesting new pattern dynamics.
\vskip 5mm


\section{Introduction}\label{sec:Introduction}
Collective behavior of large animal groups with self-organization
into robust complex patterns is a fascinating natural phenomenon \cite{Parrish-Edelstein-Keshet, Krause-Ruxton}.
Prominent examples are schools of fish, which can move in a rather ordered formation, change direction abruptly,
or, under predator threat, swirl like a vehemently stirred fluid \cite{Pitcher-Parrish}.
Flocks of hundreds of starlings can fly as a uniformly moving group,
as well as produce turbulent aerial displays \cite{Emlen}.
Apart from its biological and evolutionary relevance,
collective phenomena play a prominent role in many other scientific
disciplines, such as robotics, control theory,
economics and social sciences, see, e.g., the recent surveys \cite{Carrillo-review, Vicsek-survey}.

Regarding its pratical and theoretical importance,
it is not surprising that the topic attracted
wide attention of physical and mathematical scientific
communities. Many different types of mathematical models
have been proposed and studied during past decades.
Their common assumption is that permanently moving
self-propelled individuals (agents) interact and adapt
their behavior according to their conspecifics.
Among them, the Cucker-Smale model is very well
known, in particular due to its relative simplicity
and the implied convenience for analytical study.
In fact, the original motivation in \cite{CS1, CS2},
where the model was introduced, was to describe
language evolution. Only in subsequent studies
it has been eventually related mainly to the description
of the emergence of flocking in groups of interacting animals.
The model considers a group of $N\in\N$ agents with
time-dependent locations $x_i(t)\in\R^d$ and velocities
$v_i(t)\in\R^d$, subject to
\(
   \dot x_i &=& v_i \,,  \label{CSm1}\\
   \dot v_i &=& \frac{1}{N} \sum_{j\neq i} g_{ij} (v_j-v_i) \,,\qquad\mbox{for } i=1,\dots,N. \label{CSm2}
\)
In the classical setting the \emph{communication rates} $g_{ij}$ depend on the metric distance
between the $i$-th and $j$-th individual, $g_{ij} = g(|x_i-x_j|)$ with a fixed function
$g: [0,\infty) \to [0,\infty)$. In particular, the rate function considered in \cite{CS1, CS2}
and most of the subsequent papers is
\(   \label{g_metric}
    g(s) = \frac{\lambda}{(\sigma^2+s^2)^\beta},
\)
with $\lambda$, $\sigma$ and $\beta$ positive parameters. Then, it was shown that
if $\beta < 1/2$, the model exhibits the so-called \emph{unconditional flocking},
where for every initial configuration
the velocities $v_i(t)$ converge to the common \emph{consensus value}
$\frac1N \sum_{i=1}^N v_i(0)$ as $t\to\infty$.
On the other hand, with $\beta\geq 1/2$ the flocking is \emph{conditional},
i.e., the asymptotic behavior of the system depends on the values of $\lambda$ and $\sigma$
and on the initial configuration. 
This result was first proved in \cite{CS1, CS2} using tools from graph theory
(spectral properties of graph Laplacian), and slightly later reproved in \cite{Tadmor-Ha} by means of elementary
calculus. Another proof has been provided in \cite{Ha-Liu}, based on bounding \eqref{CSm1}--\eqref{g_metric}
by a system of dissipative differential inequalities, and, finally, the proof of \cite{CFRT}
is based on bounding the maximal velocity.
Moreover, the case of singular communication rate $g(s) = \lambda / s^\beta$ was studied in \cite{Ha-Liu}.

Another classical topic of study is the derivation and analysis of the mean-field limit as $N\to\infty$,
which leads to the Vlasov-type kinetic equation
\(   \label{Vlasov}
    \partial_t f + \grad_x\cdot (vf) + \grad_v\cdot(\calG[f]f) = 0,
\)
with
\[
    \calG[f](t,x,v) = -\int_{\R^{2d}} g(|x-y|) (v-w) f(t,y,w) \d w\d y,
\]
where $f=f(t,x,v)$ is the agent distribution function; see, e.g.,
\cite{CCR, CFRT, Ha-Liu}. A stochastic version of the model was considered in \cite{HLL}.

In our paper we introduce a modification of \eqref{CSm1}--\eqref{CSm2}
where the communication rates $g_{ij}$ depend on the topological configuration of the group,
rather than the metric distance between the agents $i$ and $j$.
This is inspired by recent extensive observations of starling flocks
in \cite{Balerrini}, which indicate that starling interactions are scale-free
and their strength depends on the so-called ``topological distance'' between individuals,
measured in units of average bird separation.
In other words, the relevant quantity is how many intermediate individuals separate
two birds, not how far apart they are in the metric sense.
Moreover, it is the very shape of the interaction that
depends on the topological distance, not simply the cut-off
or the range. Clearly, the word "distance" is a misnomer,
since the mathematical definition of a \emph{distance} is violated in two ways:
the symmetry and triangle inequality do not hold. Therefore, let us introduce
the more appropriate terminology \emph{relative separation} of the $i$-th and $j$-th agent,
defined as
\(   \label{topdist}
    \td_{ij} = \#\bigl\{1\leq k \leq N; |x_i-x_k| < |x_i-x_j| \bigr\} = \sum_{k=1}^N \chi_{[0,1)}\left(\frac{|x_i-x_k|}{|x_i-x_j|}\right),
\)
where $\chi_{[0,1)}$ is the characteristic function of the real interval $[0,1)$.
Observe that $\td_{ij}=0$ if and only if $x_i=x_j$.
We will then consider the model \eqref{CSm1}--\eqref{CSm2} with the communication rates $g_{ij}$
depending on the relative separation $\td_{ij}$ through $g_{ij} = g(\td_{ij})$,
where $g: \N \to [0,\infty)$ is a given function.
In our paper, we are interested in studying the asymptotic flocking behavior
of this model in dependence of the function $g$ and the configuration
of the agents. We will show that fundamental for the topological model
to exhibit asymptotic flocking is the connectivity of the graphs of inter-agent interactions.
Moreover, we will formally derive the mean-field limit of the model,
which will be of the Vlasov-type \eqref{Vlasov}.
Finally, we will propose that the idea of introducing
scale-free interactions is applicable to a broad spectrum of discrete and continuum
models of collective behavior, and mention the well-known attraction-repulsion model
\cite{attr-rep} as a particular example.
  
Before we proceed to the mathematical analysis of the topological Cucker-Smale
model, let us make a short remark about the biological relevance
of topological interactions, as opposed to metric ones.
Indeed, as argued in \cite{Balerrini}, the difference
between metric and topological interactions has a major
impact on the global flocking behavior of the corresponding models,
and, consequently, on their biological relevance.
The most important evolutionary advantage of collective behavior
is avoidance of predation, for which strong cohesion of the group
is essential \cite{Pitcher-Parrish}. Numerical experiments performed in \cite{Balerrini} and \cite{Bode}
strongly suggest that only topological interactions grant
such robust cohesion and, therefore, high biological fitness.
The intuitive explanation of this fact is the high ``flexibility''
of topological interactions, compared to metric ones:
In the metric case, interaction effectively vanishes and cohesion is lost
whenever the interindividual distance becomes larger then the prescribed radius.
In contrast, topological interaction stays effective
also over long distances or in the case of low individual density,
and so keeps the cohesion even in the presence of strong perturbations,
of which predation is the most relevant.
Another work revealing the qualitative difference
between the two kinds of interaction is \cite{Ginelli-Chate},
where a topological Voronoi Vicsek-type model 
exhibited an ordered moving phase with novel
long-range correlations, even though no direct
long-range connections emerged, and the transition to collective motion
exhibited critical properties different from known universality classes.

The biological relevance of topological interactions contrasts
with the fact that the mathematical flocking community so far focused
almost exclusively on models with metric interactions.
This is probably best explained by their accessibility
to mathematical analysis and studies of global quantitative properties.
In our paper we aim to make a step towards closing this gap by studying
the flocking properties of the discrete topological Cucker-Smale model
and deriving its kinetic and hydrodynamic limits.
The paper is organized as follows: In Section \ref{sec:Formulation}
we offer a detailed formulation of the model,
relate it to known models of control theory
and discuss its main properties.
We also provide three examples pointing out that the connectivity
of the graphs of inter-agent interactions play a central role
for the asymptotic flocking behavior.
In Section \ref{sec:Flocking} we provide our main result
about the unconditional flocking.
In Section \ref{sec:MeanField} we derive
the kinetic and hydrodynamic limits of the model
as the population size $N$ tends to infinity and briefly
discuss their flocking properties.
Finally, in Section \ref{sec:OtherModels} we argue that our approach is applicable
to a broad class of flocking models.
As a particular example, we mention the attraction-repulsion model \cite{attr-rep}
and illustrate with numerical simulations that the modified model produces
interesting new pattern dynamics.

\section{Formulation of the model and its main properties}\label{sec:Formulation}
We consider the Cucker-Smale-type model
\(
   \dot x_i &=& v_i,\label{CSt1}\\
   \dot v_i &=& \frac{1}{\gamma_N} \sum_{j=1}^N g(\td_{ij}) (v_j-v_i) \,,\qquad\mbox{for } i=1,\dots,N, \label{CSt2}
\)
with the relative separation $\td_{ij}$ defined by \eqref{topdist},
and $g: \N \to [0,\infty)$ is a given function.
From the modeling point of view, it makes sense to assume that interactions
with closer (in the sense of the relative separation) agents are stronger,
which corresponds to an nonincreasing $g$;
however, we do not require the validity of this assumption for the forthcoming analysis.

Note that the right-hand side in \eqref{CSt2} is scaled by $1/\gamma_N$, instead of the scaling by $1/N$ in \eqref{CSm2}, with
\(  \label{gamma_N}
   \gamma_N = \sum_{i=1}^{N} g(i-1).
\)
This is motivated by the fact that $\sum_{j=1}^N g(\td_{ij}) = \gamma_N$ for every $i=1,\dots,N$,
so that $\frac{1}{\gamma_N} \sum_{j=1}^N g(\td_{ij}) v_j$ is a convex combination of the vectors $v_j$.
The same type of scaling was considered in \cite{Motsch-Tadmor}, going under the name of \emph{relative distance};
however, it is still a \emph{metric} distance, with $g_{ij}=g(|x_i-x_j|)$.
The advantage of the model \cite{Motsch-Tadmor} is that it does not involve any explicit dependence
on the number of agents; just their geometry in phase space is taken into account.
Therefore, it cures the drawback of the classical Cucker-Smale model,
where the motion of an agent is modified by the total number of agents even if its dynamics is only influenced
by essentially a few close neighbors. Our approach, which is to replace the metric distance
by its topological counterpart, treats this drawback in a much more radical way,
and, as explained in Section \ref{sec:Introduction}, has the advantage of being supported
by actual observations of biological collective motion (starlings).

The mathematical disadvantage of introducing relative or relative separations is that the model
looses its symmetry. Indeed, observe that in the classical Cucker-Smale model \eqref{CSm1}--\eqref{CSm2}
we have $g_{ij} = g(|x_i-x_j|) = g_{ji}$. Symmetry is the cornerstone for studying the long time behavior
of its solutions, since it implies that the total momentum $V(t):=\frac1N \sum_{i=1}^N v_i(t)$ is conserved,
and with appropriate assumptions one can prove the decay of the fluctuations of the velocities about $V(t)\equiv V(0)$.
In fact, the sufficient and necessary condition for momentum conservation is that the matrix of communication rates
$G=(g_{ij})\in\R^{N\times N}$ is \emph{balanced}, i.e., its row and column sums are equal,
\[
    \sum_{i=1}^N g_{ij} = \sum_{j=1}^N g_{ij}.
\]
Neither the model of \cite{Motsch-Tadmor} nor our topological model has balanced $G$. Therefore,
new analytical techniques to study the flocking behavior need to be developed.
The approach of \cite{Motsch-Tadmor} was based on the notion of active sets
and estimation of the maximal action of antisymmetric matrices.
Our approach is based on rather simple analytic arguments
and relies on the assumption of the strong connectivity of the underlying directed graph.

For a fixed $N$, let us rescale the model such that $\gamma_N=1$.
Then \eqref{CSt2} can be written in the form
\[
   \dot v = (G-I) v = -L v \,,
\]
where the rows of $v\in\R^{N\times d}$ are composed of the vectors $v_1,\dots,v_N$,
and $G=(g_{ij})\in\R^{N\times N}$ with $g_{ij}=g(\td_{ij})$ is a stochastic matrix (in the sense that all row sums are equal to $1$).
The matrix $L := I-G$ is then the discrete Laplacian
corresponding to the \emph{directed} graph (called \emph{digraph} in the sequel) spanned by the agents with edge weights given by $g_{ij}\geq 0$.
In this sense we will indetify the topological configuration of the system with the matrix $G$ or, equivalently, $L$.
Obviously, with $N$ agents one can have at most $N!$ different topologies $G^k$, $k=1,\dots,N!$, and the system switches between
(some of) them during its temporal evolution. Because the trajectories $x_i$ of the agents are continuous,
in any time interval there can be only countably many switches and any two consecutive switches
are separated by a nonempty open time interval, where the configuration does not change.
We can then understand the system \eqref{CSt1}--\eqref{CSt2} as a switching communication network,
see, e.g., \cite{Olfati-Saber-Murray},
\(
   \dot x(t) &=& v(t) \,,  \label{SN1}\\
   \dot v(t) &=& -L^{\sigma(x(t))} v(t) \,, \label{SN2}
\)
with $-L^{\sigma(x(t))} = G^{\sigma(x(t))} - I$ and the function $\sigma: \R^{N\times d} \to \{1,\dots,N!\}$
switches between the topologies according to the current configuration $x(t)$.
We introduce the following notion of solution:

\begin{defin}\label{def:sol}
Let the curve $(x(t),v(t))_{t\geq 0} \in\R^{N\times d}\times\R^{N\times d}$
be globally continuous. Denote $\sigma: \R^{N\times d} \to \{1,\dots,N!\}$ a switching function
with at most countably many switches on $(0,\infty)$,
and denote $\{\I_k\}_{k\in\N}$ the system of open intervals
such that $\sigma\equiv\mbox{const.}$ on every $\I_k$ and $\bigcup_{k\in\N} \overline{\I_k} = [0,\infty)$.
We call the continuous curve $(x(t),v(t))_{t\geq 0}$ a solution to \eqref{SN1}--\eqref{SN2}
if it solves the differential equation on every open interval $\I_k$, $k\in\N$.
\end{defin}

For the case of a \emph{fixed} configuration $G=(g_{ij})$
with a \emph{strongly connected} digraph (let us recall that strong connectivity means
that every pair of vertices $u$, $w$ is connected by a \emph{directed} path from $u$ to $w$
and a \emph{directed} path from $w$ to $u$), Corollaries 1 and 2 of \cite{Olfati-Saber-Murray}
imply that an asymptotic consensus is found even if $G$ is not balanced:

\begin{proposition}\label{prop:OSM}
Consider the communication network $\dot v = (G-I)v$ with a fixed topology $G$
that is a strongly connected digraph. Let $\xi=(\xi_1,\dots,\xi_N)^T$ be a nonnegative left eigenvector
corresponding to the zero eigenvalue of $L=-(G-I)$, i.e., $\xi^T L = 0$,
such that $\sum_{i=1}^N \xi_i > 0$. Then, an asymptotic consensus $v^\infty\in\R^d$ is reached
exponentially fast as $t\to\infty$, i.e.,
\[
    \lim_{t\to\infty} v_i(t) = v^\infty,\qquad\mbox{for } i=1,\dots,N,
\]
with $v^\infty$ being the convex combination of $v_i$,
\(   \label{v^infty}
    v^\infty = \frac{\sum_{i=1}^N \xi_i v_i(0)}{\sum_{i=1}^N \xi_i} .
\)
\end{proposition}

The proof \cite{Olfati-Saber-Murray} is based on the simple observation
that $\xi^T v$ is an invariant quantity (note that the graph Laplacian has always a zero eigenvalue).
Let us point out the two fundamental assumptions of the above Theorem: (i) the topology is fixed and
(ii) the corresponding digraph is strongly connected.
Of course, the above result would be applicable for our system \eqref{SN1}--\eqref{SN2}
as soon as we knew that the system undergoes only a finite number of switches.
We would then simply pick the last attained configuration $G^\sigma$, and if it was strongly connected,
we would conclude that an asymptotic consensus will be reached exponentially fast as $t\to\infty$. However, the following
example shows that a special choice of the communication rate $g$ together with
a particular initial configuration leads to an infinite number of switches:

\begin{example}\label{ex1}
Let us consider a group of $7$ agents, denoted by $i=-3,-2,-1,0,1,2,3$, moving on the real line,
with initial positions
\[
    x_{-1}(0) = -6,\quad &x_{-2}(0) = -9,&\quad x_{-3}(0) = -10,\\
    x_{1}(0) = 6,\quad &x_{2}(0) = 9,&\quad x_{3}(0) = 10,
\]
and initial velocities $v_i=-1$ for $i=-1,-2,-3$ and $v_i=1$ for $i=1,2,3$.
For the agent $i=0$ we prescribe the initial datum
\(   \label{x0v0_IC}
    x_0(0) = 0,\quad v_0(0) = c,
\)
for some $c>0$.
Let us consider the model \eqref{CSt1}--\eqref{CSt2} with weight $g$ given by
\[
    g(2) = 1,\qquad g(i) = 0 \quad\mbox{for } i\neq 2,
\]
i.e., evert agents interacts exclusively with its \emph{second} closest neighbor.
Let us assume for a moment that we are able to choose $c>0$ such that $-1 \leq x_0(t) \leq 1$
for all $t\geq 0$; we will justify this assumption later.
Then, the agents from the two triplets $i\in\{-1,-2,-3\}$ and $i\in\{1,2,3\}$
will interact only with agents from the same group, so that
the velocities of these six will remain constant and equal to their initial values.
The agent $i=0$ will interact with $i=-1$ if $x_0(t)>0$ and
with $i=1$ if $x_0(t)<0$. Therefore,
\(   \label{ODEx0v0}
    \dot x_0 = v_0,\qquad \dot v_0 = -\sign(x_0) -v_0.
\)
Due to the initial condition \eqref{x0v0_IC}, the agent first moves to the right,
so that $\sign(x_0)=1$. Then \eqref{ODEx0v0} is resolved by
\(   \label{x0v0}
    x_0(t) = (c+1)(1-e^{-t})-t,\qquad v_0(t) = (c+1)e^{-t}-1,
\)
Note that since $x_0'(0) = c>0$ and $\lim_{t\to\infty} x_0(t) = -\infty$,
there exists a positive time $\tau=\tau(c)$ such that $x_0(\tau(c))=0$,
i.e., the agent will return to the origin eventually,
and then its trajectory will be be subject to \eqref{ODEx0v0}
with $\sign(x_0)=-1$.
The ``turning point'' $t_t$, when $v_0(t_t)=0$,
is $t_t = \ln(c+1)$ and $x_0(t_t) = c - \ln(c+1)$, so for $c$ small enough
the agent will be confined to the strip $-1 < x_0 < 1$, which verifies the above assumption.
The connectivity diagram of the two configurations with $x_0>0$ and $x_0<0$
is visualised in Figure \ref{fig:ex1_connectivity}.

\psset{xunit=0.6cm,yunit=0.6cm,runit=0.6cm,%
nodesep=3pt}
\def\Cnodeput(#1)#2#3{\cnode(#1){4mm}{#2}\rput(#2){#3}}

\noindent
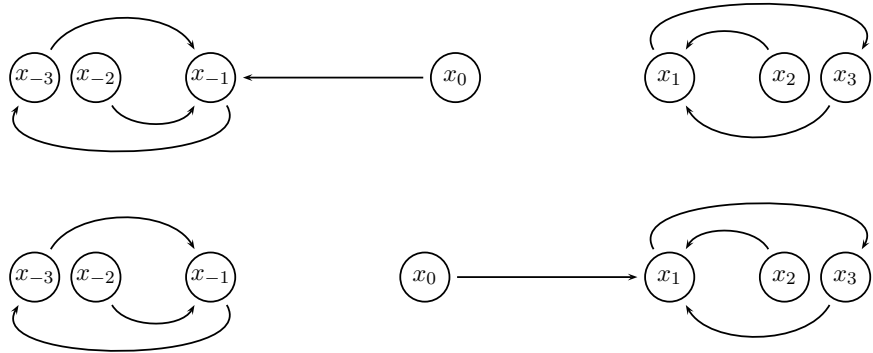
\begin{figure}
{\centering \begin{tabular}[h]{c}
\resizebox*{0.7\linewidth}{!}{
\begin{pspicture}(-12,-2)(12,2)
\Cnodeput(-10.6,0){A3}{$x_{-3}$}
\Cnodeput(-9,0){A2}{$x_{-2}$}
\Cnodeput(-6,0){A1}{$x_{-1}$}
\Cnodeput(0.4,0){A0}{$x_{0}$}
\Cnodeput(6,0){B1}{$x_{1}$}
\Cnodeput(9,0){B2}{$x_{2}$}
\Cnodeput(10.6,0){B3}{$x_{3}$}
\ncarc[arcangle=60]{->}{A3}{A1}
\ncarc[arcangle=-60]{->}{A2}{A1}
\ncarc[arcangle=120]{->}{A1}{A3}
\ncarc[arcangle=0]{->}{A0}{A1}
\ncarc[arcangle=60]{->}{B3}{B1}
\ncarc[arcangle=-60]{->}{B2}{B1}
\ncarc[arcangle=120]{->}{B1}{B3}
\end{pspicture}}
\\ $\;$\\
\resizebox*{0.7\linewidth}{!}{
\begin{pspicture}(-12,-2)(12,2)
\Cnodeput(-10.6,0){A3}{$x_{-3}$}
\Cnodeput(-9,0){A2}{$x_{-2}$}
\Cnodeput(-6,0){A1}{$x_{-1}$}
\Cnodeput(-0.4,0){A0}{$x_{0}$}
\Cnodeput(6,0){B1}{$x_{1}$}
\Cnodeput(9,0){B2}{$x_{2}$}
\Cnodeput(10.6,0){B3}{$x_{3}$}
\ncarc[arcangle=60]{->}{A3}{A1}
\ncarc[arcangle=-60]{->}{A2}{A1}
\ncarc[arcangle=120]{->}{A1}{A3}
\ncarc[arcangle=0]{->}{A0}{B1}
\ncarc[arcangle=60]{->}{B3}{B1}
\ncarc[arcangle=-60]{->}{B2}{B1}
\ncarc[arcangle=120]{->}{B1}{B3}
\end{pspicture}}
\end{tabular}\par}
\caption{Connectivity diagram for Example \ref{ex1}, upper panel for $x_0>0$ and lower panel for $x_0<0$.
An arrow is pointing from $x_i$ to $x_j$ if and only if $g_{ij}>0$.}
\label{fig:ex1_connectivity}
\end{figure}

However, the agent will never exhibit a periodic trajectory,
because $v_0(\tau(c)) > -c$.
Indeed, since $x_0(t)+v_0(t) = c-t$, we have $v_0(\tau(c)) = c-\tau(c)$.
Now if there was a $c>0$ such that $v_0(\tau(c))=-c$, this
would imply $\tau(c) = 2c$ and inserting this into \eqref{x0v0}, we would get
\[
    x_0(\tau) = (1+\tau/2)(1-e^{-\tau})-\tau,
\]
which, as can be easily checked, is equal to zero only for $\tau=0$.
Therefore, the movement of $x_0$ is never periodic,
but still the agent crosses the origin $x=0$ infinitely many times,
thus changing the configuration of the system infinitely many times.
To see this, let us calculate the implicit derivative of $x_0(\tau(c))=0$ with respect to $c$,
\[
    \tau'(c) = \frac{e^{-\tau(c)}-1}{(c+1)e^{-\tau(c)}-1} = \frac{e^{-\tau(c)}-1}{c-\tau(c)}.
\]
Realizing that $v_0(\tau(c)) = c-\tau(c) \leq 0$ (the agent can only return back to the origin
with a nonpositive velocity), we find $\tau'(c) > 0$ for all $c > 0$. 
Moreover, using the L'Hospital rule, we calculate
\[
    \tau'(0) = 2,\qquad \tau''(0) = -4/3.
\]
This implies that for $c$ small enough, $v_0(\tau(c)) = c-\tau(c) \geq -c$,
i.e., the speed $s(c):=|v_0(\tau(c))|$ of the agent when returning to the origin
is smaller than its initial speed $c$ when it was leaving the origin.
Moreover, the time needed to return to the origin, $\tau(c)$, satisfies $c < \tau(c) < 2c$ for all $c>0$.
A plot of the trajectories and velocities of the agents is provided in Fig. \ref{fig:ex1}.

\noindent
\begin{figure}
{\centering \begin{tabular}[h]{cc}
\resizebox*{0.45\linewidth}{!}{\includegraphics{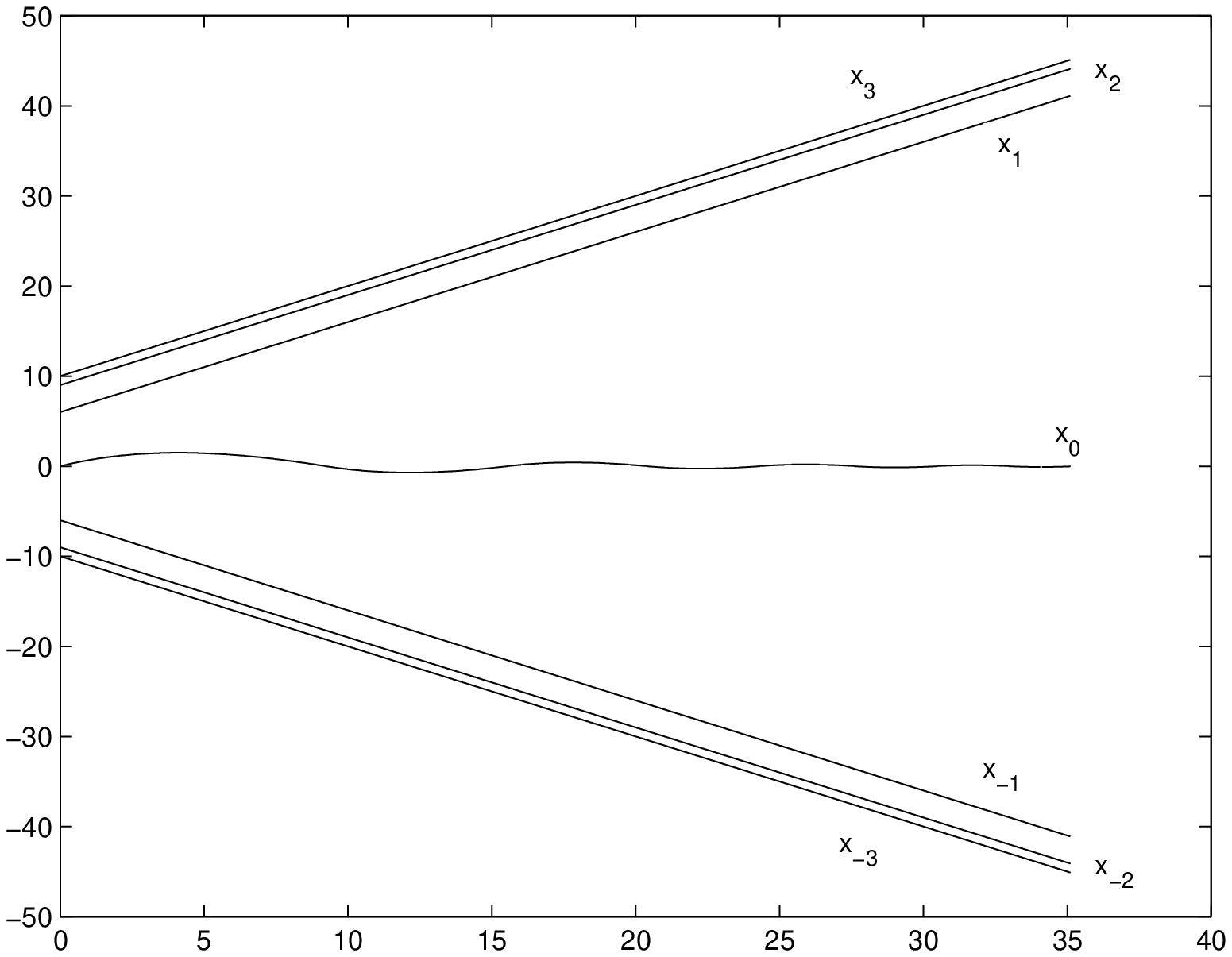}} &
\resizebox*{0.45\linewidth}{!}{\includegraphics{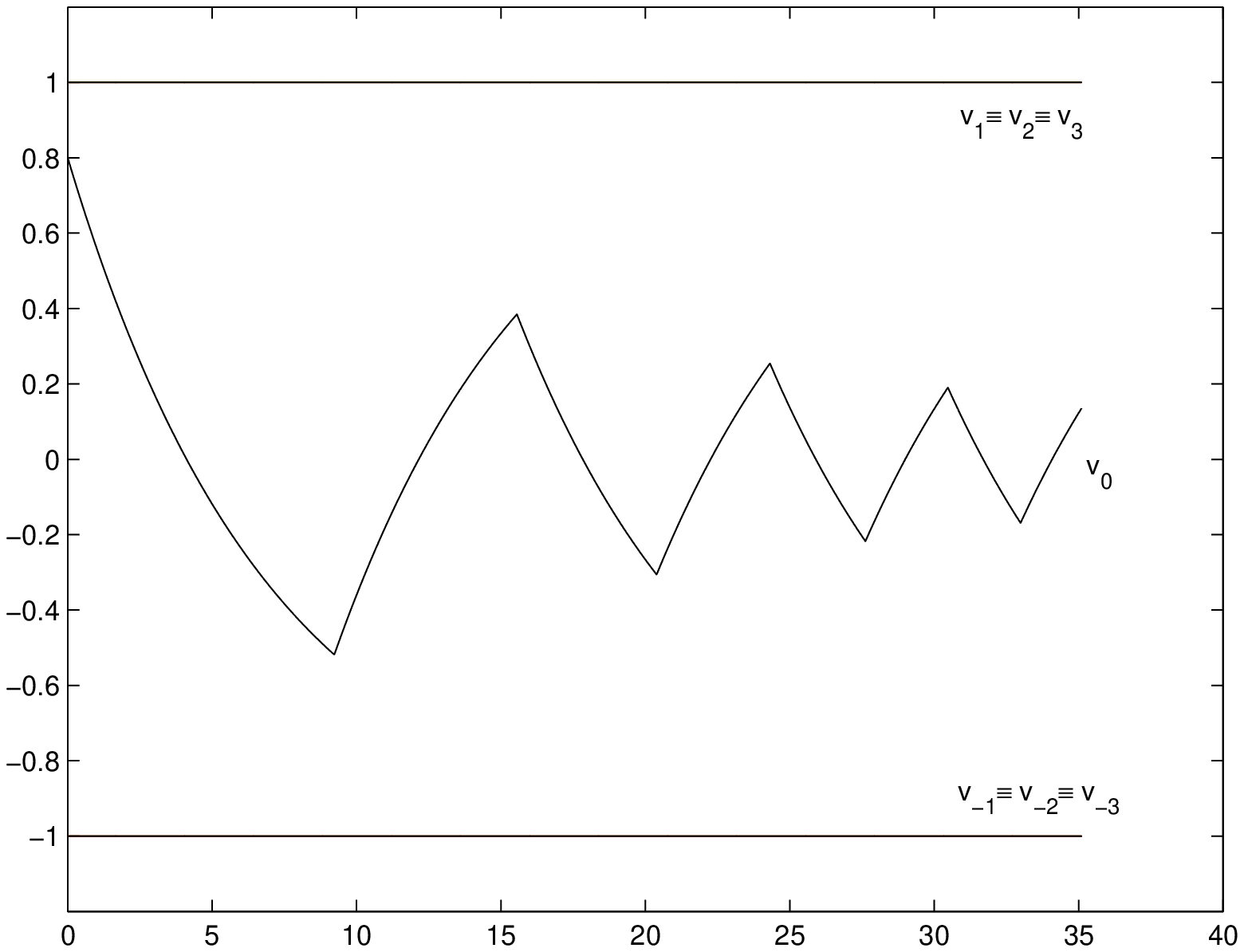}}
\end{tabular}\par}
\caption{Agent locations (left panel) and velocities (right panel)
in Example \ref{ex1}.}
\label{fig:ex1}
\end{figure}

Thus, the agent follows a nonperiodic trajectory, starting from $x_0(0)=0$ and $v_0(0)=c>0$,
and, since $s(c)>0$ for all $c>0$, returning to the origin infinitely many times,
at instances $\tau(c)$, $\tau(s(c))$, $\tau(s(s(c)))$, $\dots$, with speeds,
respectively, $s(c)$, $s(s(c))$, $s(s(s(c)))$, $\dots$.
Consequently, the configuration of the system switches infinitely many times.
However, the fact that both the sequence of return times and the sequence of speeds are strictly decreasing
(and tend to zero, since otherwise a periodic solution would exist)
immediately poses the question whether the series of return times diverges or not.
If it does diverge, then the nonperiodic solution of the ODE \eqref{ODEx0v0}
never stops to oscillate. On the other hand, a convergent series of return times
would mean that there exists a $T>0$ such that $x_0(t)\equiv v_0(t) \equiv 0$ for all $t\geq T$. 
Unfortunately, we were not able to find an analytic answer to this question;
let us remark that since
\[
    \lim_{c\to 0} \frac{s(c)}{c} = \lim_{c\to 0} \frac{\tau(c)-c}{c} = \tau'(0) - 1 = 1,
\]
the ratio test for convergence of series is inconclusive,
and we were not able to obtain the answer from any of the higher-order tests.
The only trivial observation is that the answer does not depend on the particular value of $c$.
\end{example}

Our second example shows that even if only a finite number of switches takes place,
\emph{strong} connectivity of the digraph of the last attained configuration is necessary
for a consensus to be found. Indeed, if the digraph is only \emph{weakly} connected
(let us recall that a digraph is called weakly connected if replacing all
of its directed edges with undirected edges produces a connected undirected graph),
the system may not find any consensus:

\begin{example}\label{ex2}
Let us consider a group of $7$ agents, denoted by $i=-3,-2,-1,0,1,2,3$, moving on the real line,
with initial positions
\[
    x_{-1}(0) = -1,\quad &x_{-2}(0) = -3,&\quad x_{-3}(0) = -6,\\
    x_{1}(0) = 1,\quad &x_{2}(0) = 3,&\quad x_{3}(0) = 6,
\]
and initial velocities $v_i=-1$ for $i=-1,-2,-3$ and $v_i=1$ for $i=1,2,3$.
For the agent $i=0$ we prescribe the initial datum
\[
    x_0(0) = 0,\quad v_0(0) = 0.
\]
Let us consider the model \eqref{CSt1}--\eqref{CSt2} with weight $g$ given by
\[
    g(1) = g(2) = 1/2,\qquad g(i) = 0 \quad\mbox{for } i \geq 3,
\]
i.e., every agents interacts with its two closest neighbors, updating its velocity
according to the arithmetic mean of their velocities.

\noindent
\begin{figure}
{\centering \begin{tabular}[h]{cc}
\resizebox*{0.49\linewidth}{!}{\includegraphics{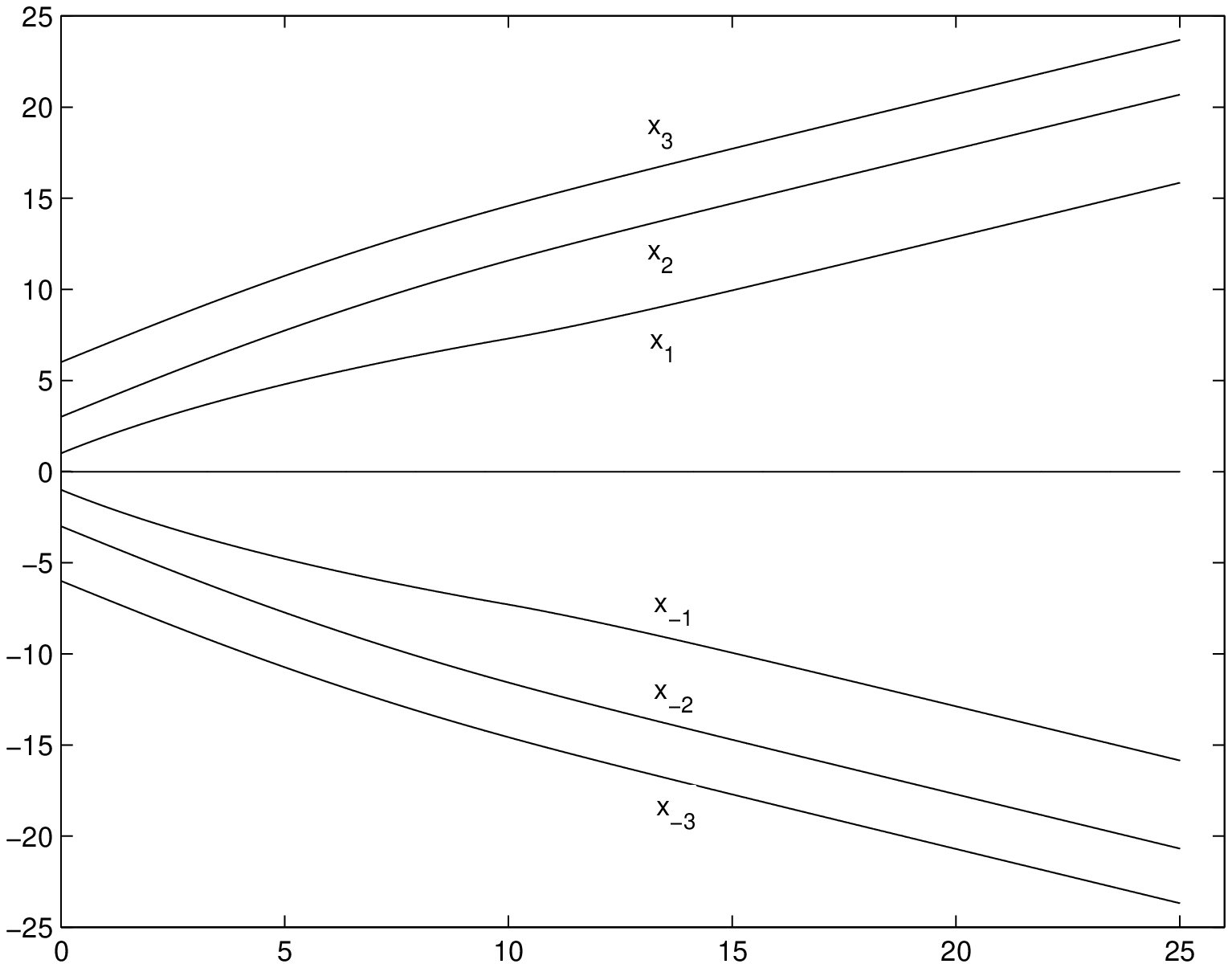}} &
\resizebox*{0.49\linewidth}{!}{\includegraphics{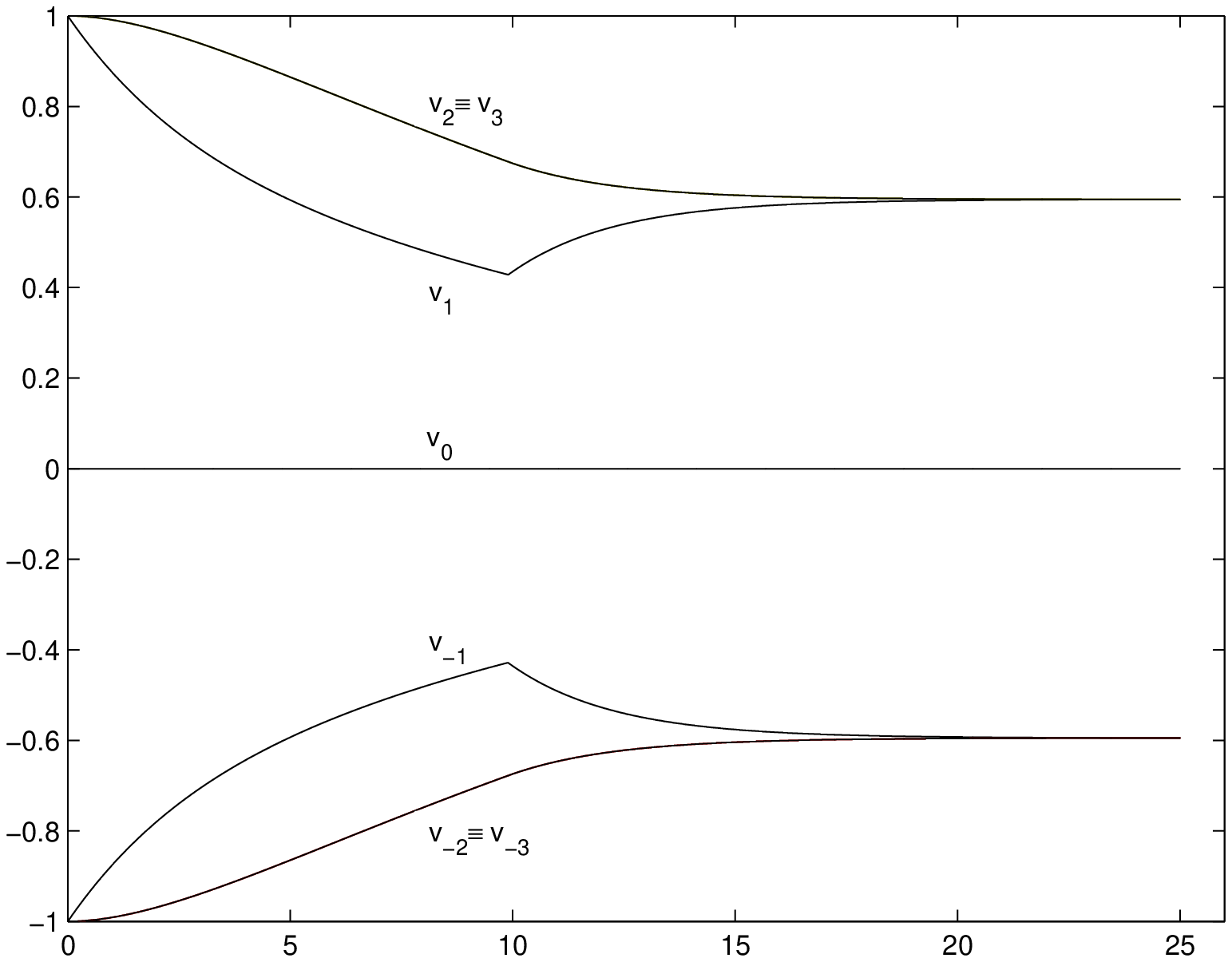}}
\end{tabular}\par}
\caption{Agent locations (left panel) and velocities (right panel)
in Example \ref{ex2}.}
\label{fig:ex2}
\end{figure}

Instead of attempting to find an analytic solution to the corresponding ODE system,
we performed a numerical simulation (see Fig. \ref{fig:ex2}), showing that the system takes two configurations:
\begin{itemize}
\item First, for $t \leq t_0$, with $t_0 \simeq 10$, the system stays in the initial configuration
with a strongly connected digraph, plotted in the upper panel of Fig. \ref{fig:ex2_connectivity}.
The configuration remains symmetric with $x_0 \equiv 0$, $v_0 \equiv 0$.
\item After time $t_0$, the distance between $x_1$ and $x_3$ becomes smaller than
the distace between $x_1$ and $x_0$. Due to the symmetry,
the same happens for $x_{-1}$ and $x_{-3}$. Therefore, the connectivity of the digraph
changes from strong to weak, see lower panel of Fig. \ref{fig:ex2_connectivity}.
Then, the agents of the group $\{1,2,3\}$ interact only among themselves,
and the same holds for the group $\{-1,-2,-3\}$.
Therefore, each of the two groups will find their own velocity consensus as $t\to\infty$,
which will be nonzero velocities with opposite signs.
The zeroth agent will stay with $x_0 \equiv 0$, $v_0 \equiv 0$.
Consequently, no global velocity consensus will be achieved.
\end{itemize}
Let us note that in this example the initial digraph was strongly connected,
but this property got lost during the evolution, due to switching to another configuration.
This new configuration is only weakly connected, which is not sufficient
for a velocity consensus to be found.

\psset{xunit=1cm,yunit=1cm,runit=1cm,%
nodesep=3pt}
\def\Cnodeput(#1)#2#3{\cnode(#1){4mm}{#2}\rput(#2){#3}}

\noindent
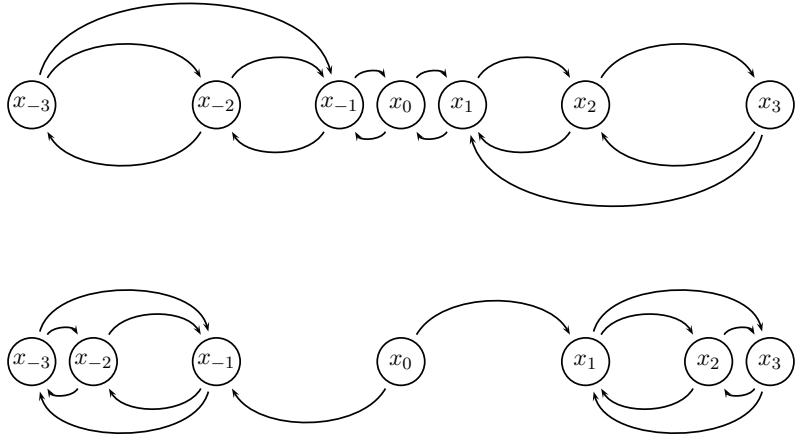
\begin{figure}
{\centering \begin{tabular}[h]{c}
\resizebox*{0.75\linewidth}{!}{
\begin{pspicture}(-8,-2)(8,2)
\Cnodeput(-6,0){A3}{$x_{-3}$}
\Cnodeput(-3,0){A2}{$x_{-2}$}
\Cnodeput(-1,0){A1}{$x_{-1}$}
\Cnodeput(0,0){A0}{$x_{0}$}
\Cnodeput(1,0){B1}{$x_{1}$}
\Cnodeput(3,0){B2}{$x_{2}$}
\Cnodeput(6,0){B3}{$x_{3}$}
\ncarc[arcangle=60]{->}{A0}{B1}
\ncarc[arcangle=60]{->}{B1}{B2}
\ncarc[arcangle=60]{->}{B2}{B3}
\ncarc[arcangle=60]{->}{B3}{B2}
\ncarc[arcangle=60]{->}{B2}{B1}
\ncarc[arcangle=75]{->}{B3}{B1}
\ncarc[arcangle=60]{->}{B1}{A0}
\ncarc[arcangle=60]{->}{A0}{A1}
\ncarc[arcangle=60]{->}{A1}{A2}
\ncarc[arcangle=60]{->}{A2}{A3}
\ncarc[arcangle=60]{->}{A3}{A2}
\ncarc[arcangle=60]{->}{A2}{A1}
\ncarc[arcangle=75]{->}{A3}{A1}
\ncarc[arcangle=60]{->}{A1}{A0}
\end{pspicture}}\\
\resizebox*{0.75\linewidth}{!}{
\begin{pspicture}(-8,-2)(8,2)
\Cnodeput(-6,0){A3}{$x_{-3}$}
\Cnodeput(-5,0){A2}{$x_{-2}$}
\Cnodeput(-3,0){A1}{$x_{-1}$}
\Cnodeput(0,0){A0}{$x_{0}$}
\Cnodeput(3,0){B1}{$x_{1}$}
\Cnodeput(5,0){B2}{$x_{2}$}
\Cnodeput(6,0){B3}{$x_{3}$}
\ncarc[arcangle=60]{->}{A0}{B1}
\ncarc[arcangle=60]{->}{B1}{B2}
\ncarc[arcangle=60]{->}{B2}{B3}
\ncarc[arcangle=60]{->}{B3}{B2}
\ncarc[arcangle=60]{->}{B2}{B1}
\ncarc[arcangle=75]{->}{B3}{B1}
\ncarc[arcangle=75]{->}{B1}{B3}
\ncarc[arcangle=60]{->}{A0}{A1}
\ncarc[arcangle=60]{->}{A1}{A2}
\ncarc[arcangle=60]{->}{A2}{A3}
\ncarc[arcangle=60]{->}{A3}{A2}
\ncarc[arcangle=60]{->}{A2}{A1}
\ncarc[arcangle=75]{->}{A3}{A1}
\ncarc[arcangle=75]{->}{A1}{A3}
\end{pspicture}}
\end{tabular}\par}
\caption{Connectivity diagram for Example \ref{ex2}. Upper panel shows the strongly connected configuration
until $t\simeq 10$, which then changes to the weakly connected configuration in the lower panel.
An arrow is pointing from $x_i$ to $x_j$ if and only if $g_{ij}>0$.}
\label{fig:ex2_connectivity}
\end{figure}

\end{example}

Finally, our third and last example shows that there is no general assumption
on the function $g$ that would guarantee strong connectivity of all the graphs
corresponding to $G=(g(\td_{ij}))$. In particular, even if $g$ is positive for all but the farthest agent,
a configuration can be easily found which is only weakly connected.

\begin{example}\label{ex3}
Let us fix $N\in\N$ and let
\[
    g(i)=1 \mbox{ for } i=1,\dots,N-2, \qquad g(N-1)=0.
\]
Let us consider a group of $N$ agents $x_1,\dots, x_N$, with
\[
    |x_i|\leq 1 \mbox{ for } i=1,\dots,N-1,\qquad |x_N|=5.
\]
Then $\td_{iN}=N-1$, and, consequently, $g(\td_{iN})=0$
for all $i=1,\dots,N-1$. Clearly, there is no oriented path
from, say, $x_1$, to $x_N$, and so the corresponding digraph
is not strongly connected.
\end{example}

\section{Asymptotic flocking}\label{sec:Flocking}
We now present our main result, which applies to the full model \eqref{CSt1}--\eqref{CSt2},
where we assume global existence of solutions. In fact, in the proof we do not explicitly make use
of the fact that the function $\sigma(x)$ switches the topological configurations $G^k$
according to the relative separations $g_{ij}=g(\td_{ij})$.
The only relevant property is that there is at most a countable number of switches,
which is guaranteed by the continuity of the trajectories $x(t)$. The fundamental assumption is then
that there exists a topological configuration with a strongly connected digraph where the system spends
an infinite amount of time. This of course is a strong assumption and might be difficult
to verify practically, however, the above Examples \ref{ex1}, \ref{ex2} and \ref{ex3} show that it is necessary.
From this point of view, the forthcoming Theorem is optimal.

\begin{theorem}\label{thm:flocking}
Let $(x(t),v(t))_{t\geq 0}$ denote the solution of the system \eqref{CSt1}--\eqref{CSt2}
in the sense of Definition \ref{def:sol},
subject to the initial condition $(x(0), v(0)) \in \R^{N\times d}\times\R^{N\times d}$.
Assume that there exists a topological configuration, say $G^0 = (g^0_{ij})$,
with a strongly connected digraph, where the system spends an infinite amount of time, i.e.,
\[
    \bigl|\{t\geq 0;\; g(\td_{ij}) = g^0_{ij} \mbox{ for all }i,j=1,\dots,N\}\bigr| = +\infty.  
\]
Then the system finds an asymptotic velocity consensus, i.e., there exists a vector $v^\infty$
in the convex hull of $\{v_1(0),\dots,v_N(0)\}$ such that
\(   \label{thm1:statement}
    \lim_{t\to\infty} v_i(t) = v^\infty\,,\qquad \mbox{for all }  i=1,\dots,N\,.
\)
\end{theorem}

\begproof
The proof will be carried out in several steps.

\textbf{Step 1: Maximal speed.}
Refering to Definition \ref{def:sol},
let us denote by $\I_k := (t_{k-1}, t_k)$, $k\in\N$, the system of disjoint open time intervals
where the system does not change its topological configuration,
and such that $\bigcup_{k\in\N} \overline{\I_k} = [0,\infty)$.
Inspired by \cite{CFRT}, let us define the function
\[
    \omega(t) := \max_{i=1,\dots,N} |v_i(t)| \qquad\mbox{for all } t\geq 0 \,,
\]
and, moreover, denote $M(t) := \mbox{argmax}_{i=1,\dots,N} |v_i(t)|$; if $M(t)$ is not uniquely determined
(i.e., there are several vectors $v_i(t)$ with maximal length), we choose one of the indices arbitrarily,
but in such a way that $M(t)$ stays constant on the longest time interval.
Since the number of particles $N$ is finite and their trajectories in the $(x,v)$-space are continuous,
there exists an at most countable system of open disjoint intervals $(\K_k)_{k\in\N}$, such that $\bigcup_{k\in\N} \overline{\K_k} = [0,\infty)$,
and $M(t)$ is constant on every $\K_k$. 
To ease the notation, we will usually skip the explicit dependence of $M$ on $t$ (or $k$) in the sequel.

By intertwining the two systems $(\I_k)_{k\in\N}$ and $(\K_k)_{k\in\N}$, we construct another at most countable system of disjoint intervals,
denoted by abuse of notation again by $(\I_k)_{k\in\N}$, such that the topological configuration does not change
\emph{and} the index $M$ is constant on each $\I_k$.
Then, on every $\I_k$ the trajectories of the particles in phase space are smooth and we can write
\[
    \frac12 \tot{}{t} \omega(t)^2 = \frac12 \tot{}{t} |v_M|^2 &=& \sum_{j=1}^N g_{Mj} (v_j - v_M)\cdot v_M \\
	  &\leq& \sum_{j=1}^N g_{Mj} \left( |v_j| - |v_M| \right) |v_M| \,. 
\]
Dividing by $|v_M|$ (note that if $|v_M|$ was zero, there would be nothing to prove), we obtain
\[
    \tot{}{t} \omega(t) = \tot{}{t} |v_M| \leq \sum_{j=1}^N g_{Mj} \left( |v_j| - |v_M| \right) \leq 0 \qquad\mbox{on every }\I_k,
\]
where the nonpositivity of the right-hand side is due to the maximality of $|v_M|$.
Observe that the above inequality holds universally, regardless of whether the configuration
is (strongly) connected or not.
Consequently, $\omega(t)$ is a globally continuous, nonincreasing and nonnegative function,
so that there exists an $0 \leq \omega_\infty \leq \omega(0)$ such that
$\lim_{t\to\infty} \omega(t) = \omega_\infty$.

\textbf{Step 2: One velocity.}
Since, by assumption, the system spends an infinite amount of time in the strongly connected configuration $G^0$,
we can pick the corresponding subsystem out of $(\I_k)_{k\in\N}$, with infinite length. By a further
subselection we get the system $\sIkn$ of infinite length, where $M(t) \equiv M_0$ for some fixed $1\leq M_0 \leq N$.
Therefore, denoting $\I^0 := \bigcup_{n\in\N}\Ikn$, we have the configuration $G^0$ and the maximal vector index $M_0$ for all $t\in\I^0$.
Moreover,
\[
    |v_{M_0}(t)| \to \omega_\infty\qquad\mbox{and}\qquad \tot{}{t} |v_{M_0}(t)|^2 \to 0 \qquad\mbox{as } t\to\infty,\; t\in\I^0,
\]
where the convergence of the time derivative is due to the monotonicity of $|v_{M_0}(t)|$.
Since the configuration $G^0$ is fixed and strongly connected,
there exists an index $j_0$ such that $g_{M_0,j_0} =: g^0 >0$ on $\I^0$.
Then we have, for $t\in\I^0$,
\[
    \frac12 \tot{}{t} |v_{M_0}(t)|^2 = \sum_{j=1}^N g_{M_0,j} (v_j-v_{M_0})\cdot v_{M_0}(t) \leq g^0 (v_{j_0}-v_{M_0})\cdot v_{M_0}(t) \leq 0,
\]
where we used the inequality $(v_j-v_{M_0})\cdot v_{M_0} \leq 0$ implied by the maximality of $|v_{M_0}|$.
Now, since the left-hand side tends to zero as $t\to\infty$ and $g^0$ is a constant, we have
\[
    v_{j_0}\cdot v_{M_0}(t) - |v_{M_0}(t)|^2 \to 0 \qquad\mbox{as } t\to\infty,\; t\in\I^0,
\]
and this futher implies $v_{j_0}\cdot v_{M_0}(t) \to \omega_\infty^2$.
Finally, since $v_{j_0}\cdot v_{M_0} = \omega_\infty^2$ if and only if $v_{j_0} = v_{M_0}$
(equality in the Cauchy-Schwartz inequality), we have
\(   \label{convergence one velocity}
    (v_{j_0} - v_{M_0})(t) \to 0 \qquad\mbox{as } t\to\infty,\; t\in\I^0.
\)
Moreover, we calculate
\[
    \frac12 \tot{}{t} |v_{j_0}(t)|^2 &=& \sum_{l=1}^N g_{j_0,l} (v_l-v_{j_0})\cdot v_{j_0}(t) \\
           &=& \sum_{l=1}^N g_{j_0,l} \left[ (v_l-v_{M_0})\cdot v_{M_0} + (v_l-v_{M_0})\cdot(v_{j_0}-v_{M_0}) + (v_{M_0}-v_{j_0})\cdot v_{j_0} \right] \\
           &\leq& \sum_{l=1}^N g_{j_0,l} (v_l-v_{M_0})\cdot v_{M_0}(t) + 3\omega(0)|v_{j_0}-v_{M_0}| \\
           &\leq& 3\omega(0)|v_{j_0}-v_{M_0}|,
\]
where we used the estimate $\max\left(|v_{j_0}|,|v_l|\right) \leq \omega(0)$ and the maximality of $v_{M_0}$.
By \eqref{convergence one velocity} we have then
\[
   \limsup_{t\to\infty,\,t\in\I^0} \tot{}{t} |v_{j_0}(t)|^2 \leq 0. 
\]
Since $|v_{j_0}(t)| \to \omega_\infty$ from below on $\I^0$ as $t\to\infty$, we conclude that
\(   \label{convergence derivative}
     \lim_{t\to\infty,\,t\in\I^1} \tot{}{t} |v_{j_0}(t)|^2 = 0.
\)
where $\I^1$ is a system of subintervals of $\I^0$, still of infinite Lebesgue measure.

\textbf{Step 3: All velocities.}
We will show that \eqref{convergence one velocity} holds for any index $\hat\jmath\in\{1,\dots,N\}$.
Due to the simple connectivity of the digraph $G^0$, there exists a path
$M_0 \mapsto j_0 \mapsto j_1 \mapsto \dots \mapsto j_\ell \mapsto \hat\jmath$, such that
$g_{M_0,j_0}>0$, $g_{j_0,j_1}>0$, $\dots$, $g_{j_\ell, \hat\jmath}>0$ on $\I^0$.
We proceed inductively, showing first that the results \eqref{convergence one velocity} and \eqref{convergence derivative} hold for $j_1$ as well.
Indeed, passing to the limit in 
\[
    \frac12 \tot{}{t} |v_{j_0}(t)|^2 &\leq& \sum_{l=1}^N g_{j_0,l} (v_l-v_{M_0})\cdot v_{M_0}(t) + 3\omega(0)|v_{j_0}-v_{M_0}| \\ 
       &\leq& g_{j_0,j_1} (v_{j_1}-v_{M_0})\cdot v_{M_0}(t) + 3\omega(0)|v_{j_0}-v_{M_0}|,
\]
we obtain, due to \eqref{convergence derivative} and the maximality of $v_{M_0}$,
\[
    0 \leq \lim_{t\to\infty,\,t\in\I^1} (v_{j_1}-v_{M_0})\cdot v_{M_0}(t) \leq 0,
\]
which immediately gives \eqref{convergence one velocity} for $v_{j_1}$ on $\I^1$.
Using this result, we argue as before to conclude
\[
   \lim_{t\to\infty,\,t\in\I^2} \tot{}{t} |v_{j_1}(t)|^2 = 0
\]
with $\I^2$ a system of subintervals of $\I^1$ of infinite Lebesgue measure.
This is \eqref{convergence derivative} for $v_{j_1}$ on $\I^2$.
Proceeding inductively, after a finite number of steps we reach the index $\hat\jmath$.

We conclude that there exists a sequence $(t_k)_{k\in\N}\subset \I^0$, $t_k\to\infty$, such that, for all $j=1,\dots,N$,
\(  \label{almost_there}
    v_j(t_k) - v_{M_0}(t_k) \to 0 \qquad\mbox{as } k\to\infty. 
\)

\textbf{Step 4: Conclusion.}
The fact that $\sum_{j=1}^N g_{ij} v_j$ is a convex combination of the velocity vectors $v_1,\dots,v_N$
(remember the scaling $\gamma_N=1$) directly implies that the convex hull of $\{v_1,\dots,v_N\}$ is nonexpanding in time,
\[
   \mbox{ch}\{v_1,\dots,v_N\}(t) \subseteq \mbox{ch}\{v_1,\dots,v_N\}(s) \qquad\mbox{for all } t>s\geq 0.
\]
Due to \eqref{almost_there}, its diameter shrinks to zero as $t\to\infty$,
and, consequently, there exists a vector $v^\infty\in\mbox{ch}\{v_1,\dots,v_N\}(0)$ such that
\[
    \lim_{t\to\infty} v_j(t) = v^\infty \qquad\mbox{for all } j=1,\dots,N. 
\]
\endproof

\begin{remark}\label{rem:position fluctuations}
The ``classical'' definitions of time-asymptotic flocking in the Cucker-Smale model,
see, e.g., \cite{CS1, CS2, Ha-Liu}, pose, in addition to the velocity alignment,
the requirement of uniform boundedness of the position fluctuations in time (formation of a group),
\(   \label{position fluctuations}
    \sup_{t\geq 0} \sum_{i=1}^N |x_i(t)-x_c(t)|^2 < \infty,
\)
where $x_c(t) = \frac1N \sum_{i=1}^N x_i(t)$ is the centre of gravity.
Such a result can be easily obtained for our topological model
under the assumption of a finite number of switches with a strongly connected
final topological configuration. Indeed, Proposition \ref{prop:OSM} provides then an
exponential convergence of the velocities to a consensus and \eqref{position fluctuations}
follows straightforwardly. However, in the general setting with possibly
infinitely many switches, we are only able to provide the estimate
\[
    \max_{i=1,\dots,N} |x_i(t)| \leq \max_{i=1,\dots,N} |x_i(0)| + \omega(0)t,
\]
which is a trivial consequence of the bound on the maximal velocity $|v_M|\leq \omega(0)$.
\end{remark}

\begin{remark}\label{rem:v_infty}
In the classical Cucker-Smale model, or any of its modifications with a balanced communication matrix,
the mean velocity $V(t) = \frac1N \sum_{i=1}^N v_i(t)$ is invariant, so that whenever asymptotic flocking takes place,
the consensus velocity is a priori given by $v^\infty = V(0)$.
In contrast to that, our model \eqref{CSt1}--\eqref{CSt2} does not seem to posses any invariants
that would allow us to predict $v^\infty$ from the initial datum, beyond the trivial fact that
$v^\infty$ is a convex combination of $v_i(0)$, $i=1,\dots,N$.
We can therefore, similarly as in \cite{Motsch-Tadmor}, consider $v^\infty$ as an \emph{emergent property} of our model,
in the sense that the asymptotic consensus $v^\infty$ is encoded in the
dynamics of the system and not just as an invariant of its initial configuration.
\end{remark}

\section{Kinetic and hydrodynamic limits}\label{sec:MeanField}
In this section we derive the mean-field limit of \eqref{CSt1}--\eqref{CSt2}
as $N\to\infty$. For this, appropriate scaling of the relative separation \eqref{topdist} is necessary,
in particular, we introduce the \emph{normalized relative separation}
\(   \label{topdist scaled}
   \td_{ij} = \frac1N \#\bigl\{1\leq k \leq N; |x_i-x_k| < |x_i-x_j| \bigr\} = \frac1N \sum_{k=1}^N \chi_{[0,1)}\left(\frac{x_i-x_k}{|x_i-x_j|}\right),
\)
where $\chi_{[0,1)}$ is the characteristic function of the interval $[0,1)$.
To derive the formal mean-field limit, we consider the system \eqref{CSt1}--\eqref{CSt2}
with a prescribed (time dependent) matrix $G=(g_{ij})$ of communication rates $g_{ij}$
with all row sums $\sum_{j=1}^N g_{ij}$ equal to
$\gamma_N:=\sum_{i=1}^{N} g((i-1)/N)$.
Introducing the empirical measure
\(   \label{empirical measure}
    f_t^N(x,v) := \frac1N \sum_{i=1}^N \delta(x-x_i(t))\delta(v-v_i(t))
\)
for $(x,v)\in\R^d\times\R^d$, and assuming that there exists a limit $f_t^N \to f$ as $N\to\infty$,
we obtain the mean-field Vlasov equation (see, e.g., \cite{Tadmor-Ha}) for $f=f(t,x,v)$,
\(   \label{kinetic}
    \part{f}{t} + v\cdot\grad_x f = \grad_v\cdot (\calG[f]f),
\)
with
\[
   \calG[f](x,v) = \frac1\gamma \int_{\R^d}\int_{\R^d} g(x,y) (w-v) f(y,w)\d w\d y \,,
\]
where $\gamma = \lim_{N\to\infty} \frac{\gamma_N}{N} = \int_0^1 g(s)\d s$.
The term $g(x,y)$ is the formal limit of $g_{ij}=g(\td_{ij})$ and is evaluated
by writing the normalized relative separation \eqref{topdist scaled} as
\[
    \td_{ij} = \int_{\R^d}\int_{\R^d} \chi_{[0,1)}\left(\frac{z-x_i}{|x_j-x_i|}\right) f^N(z,v)\d v\d z,
\]
The formal limit as $N\to\infty$ gives then the \emph{continuum relative separation} between the locations $x$ and $y$ in $\R^d$,
\(  \label{topdist-cont}
   \td[f](x,y) = \int_{\R^d}\int_{\R^d} \chi_{[0,1)}\left(\frac{|z-x|}{|y-x|}\right) f(z,v)\d v\d z.
\)
Finally, substituting $g(\td[f](x,y))$ for $g(x,y)$, we obtain
\(  \label{calG}
    \calG[f](x,v) = \frac{1}{\gamma} \int_{\R^d}\int_{\R^d} g\bigl(\td[f](x,y)\bigr) (w-v) f(y,w)\d w\d y \,.
\)

Well-posedness and stability of measure solutions to \eqref{kinetic}--\eqref{calG} can be obtained
using the tools developed in \cite{CCR} and \cite{CCcorr}, as soon as we consider a smoothened version of the continuum relative separation.
First, let us observe that \eqref{topdist-cont} can be equivalently written as
\[
   \td[f](x,y) = \int_{\R^d}\int_{\R^d} \chi_{[0,\infty)}\left(|y-x|-|z-x|\right) f(z,v)\d v\d z,
\]
where $\chi_{[0,\infty)}$ is the characteristic function of the interval $[0,\infty)$.
Then, we introduce a smoothened version $\psi$ of $\chi_{[0,\infty)}$;
in fact, we require that $\psi$ be globally Lipschitz continuous and $0\leq\psi\leq 1$.
We obtain a ``smoothened'' continuum relative separation,
\(  \label{smoothened_td}
   \tilde\td[f](x,y) = \int_{\R^d}\int_{\R^d} \psi\left(|y-x|-|z-x|\right) f(z,v)\d v\d z,
\)
and by $\tilde\calG[f]$ we denote \eqref{calG} with $\tilde\td[f]$ in place of $\td[f]$.
In order to obtain the well-posedness result, we only need to verify that $\tilde\calG[f]$
satisfies the Hypotheses 4.9 of \cite{CCR}. We will stick to the notation
$\mathcal{A}:=C([0,T]; \mathcal{P}_c(\R^d\times\R^d))$ endowed with the $1$-Wasserstein distance $\calW_1$,
where $\mathcal{P}_c(\R^d\times\R^d)$ denotes the space of probability measures with compact support in $\R^d\times\R^d$.

\begin{lemma}
\label{lem:hyp4.9}
Let $\tilde\chi$ and $g$ be Lipschitz continuous functions.
Take any $R_0>0$ and $f,h\in\mathcal{A}$ such that $\supp(f_t)\cup\supp(g_t)\subseteq B_{R_0}$ for all $t\in[0,T]$.
Then for any ball $B_R\subseteq\R^d\times\R^d$, there exists a constant $C=C(R,R_0)$ such that
\(
   \max_{t\in[0,T]} \mathrm{Lip}_R(\tilde\calG[f]) &\leq& C,\\    \label{hyp1}
   \max_{t\in[0,T]} \Norm{\tilde\calG[f]-\tilde\calG[h]}_{L^\infty(B_R)} &\leq& C \calW_1(f,h).  \label{hyp2}
\)
Here $\mathrm{Lip}_R(\tilde\calG[f])$ denotes the Lipschitz constant of $\tilde\calG[f]$ in the ball $B_R\subset\R^d\times\R^d$.
\end{lemma}

\begproof
The uniform Lipschitz continuity of $\tilde\calG[f]$ with respect to the $v$-variable is obvious.
Thereforem, for a fixed $v\in\R^d$, we estimate
\[
    \left| \tilde\calG[f](x_1,v) - \tilde\calG[f](x_2,v) \right| \leq
       \frac{\Norm{g}_{\mathrm{Lip}}}{\gamma}  \int_{\R^d}\int_{\R^d} \left| \tilde\td[f](x_1,y) - \tilde\td[f](x_2,y) \right| |v-w| f(y,w)\d w\d y.
\]
Now, since
\[
   \bigl| \tilde\td[f](x_1,y) - \tilde\td[f](x_2,y) \bigr| &\leq&
     \int_{\R^d}\int_{\R^d} \bigl| \psi(|y-x_1|-|z-x_1|) - \psi(|y-x_2|-|z-x_2|) \bigr| f(z,v) \d v\d z \\
     &\leq&
     \Norm{\psi}_{\mathrm{Lip}} \int_{\R^d}\int_{\R^d} \bigl| |y-x_1|-|z-x_1| - |y-x_2| + |z-x_2| \bigr| f(z,v) \d v\d z \\
     &\leq&
     2 \Norm{\psi}_{\mathrm{Lip}} |x_1-x_2|
\]
and due to $|v-w|\leq |v|+|w|\leq C(R,R_0)$, we obtain 
\[
    \left| \tilde\calG[f](x_1,v) - \tilde\calG[f](x_2,v) \right| \leq 2 C(R,R_0) \frac{\Norm{g}_{\mathrm{Lip}}\Norm{\psi}_{\mathrm{Lip}}}{\gamma}\, |x_1-x_2|,
\]
which directly implies \eqref{hyp1}.

For \eqref{hyp2}, let $\pi$ be an optimal transportation plan between the measures $f$ and $h$.
Then, for any $(x,v)\in B_R$, we write
\[
   \bigl| \tilde\calG[f](x,v) - \tilde\calG[h](x,v) \bigr| &=& \left| \int_{\R^{4d}} \bigl(
        g(\td[f](x,y_1))(w_1-v) - g(\td[h](x,y_2))(w_2-v) \bigr) \d\pi(y_1,w_1,y_2,w_2)  \right| \\
      &\leq&
      \Norm{g}_{\mathrm{Lip}} \int_{\R^{4d}} \bigl|\td[f](x,y_1) - \td[h](x,y_2)\bigr| |w_1-v| \d\pi(y_1,w_1,y_2,w_2) \\
      && + \Norm{g}_{L^\infty} \int_{\R^{4d}} \bigl| |w_1-v| - |w_2-v| \bigr| \d\pi(y_1,w_1,y_2,w_2) .
\]
Using again the fact the $\pi$ has marginals $f$ and $g$, we estimate
\[
   \bigl|\td[f](x,y_1) - \td[h](x,y_2)\bigr| &=& \left| \int_{\R^{4d}} \bigl(
      \psi(|y_1-x|-|z-x|) - \psi(|y_2-x|-|\xi-x|) \bigr) \d\pi(z,w_1,\xi,w_2) \right| \\
      &\leq&
      \Norm{\psi}_{\mathrm{Lip}} \int_{\R^{4d}} \bigl(
       |y_1-y_2|+|z-\xi|  \bigr) \d\pi(z,w_1,\xi,w_2) \\
      &\leq&
      \Norm{\psi}_{\mathrm{Lip}} \bigl( |y_1-y_2| + \calW_1(f,h) \bigr).
\]
Finally, since
\[
    \int_{\R^{4d}} \bigl| |w_1-v| - |w_2-v| \bigr| \d\pi(y_1,w_1,y_2,w_2) \leq \int_{\R^{4d}} |w_1-w_2|  \d\pi(y_1,w_1,y_2,w_2)
      \leq \calW_1(f,h),
\]
and $|w_1-v|\leq |w_1|+|v|\leq C(R,R_0)$, we obtain
\[
    \bigl| \tilde\calG[f](x,v) - \tilde\calG[h](x,v) \bigr| &\leq&
      \Norm{g}_{\mathrm{Lip}} \Norm{\psi}_{\mathrm{Lip}} C(R,R_0) \int_{\R^{4d}} \bigl( |y_1-y_2| + \calW_1(f,h) \bigr) \d\pi(y_1,w_1,y_2,w_2)
        + \Norm{g}_{L^\infty} \calW_1(f,h) \\
      &\leq& \bigl( 2 C(R,R_0) \Norm{g}_{\mathrm{Lip}} \Norm{\psi}_{\mathrm{Lip}} + \Norm{g}_{L^\infty} \bigr) \calW_1(f,h).
\]
\endproof

Based on the above Lemma, the theory developed in \cite{CCR} provides existence of compactly supported, global measure solutions
to the kinetic equation \eqref{kinetic} with \eqref{smoothened_td}, subject to a compactly supported measure initial condition.
These solutions are understood in the sense of push-forward measure and are unique in their class;
see Theorem 4.11 of \cite{CCR} for details.
Moreover, the solutions are stable in the $\mathcal{W}_1$-Wasserstein topology.
Since the empirical measure $f$ given by \eqref{empirical measure}, with $(x_i,v_i)$ being a solution
of the discrete topological Cucker-Smale system with $g_{ij}:=g(\tilde\td(x_i,x_j))$,
is a measure solution to the kinetic equation, the stability result justifies the rigorous
passage to the mean-field limit $N\to\infty$; see Lemma 5.1 and Corollary 5.2 of \cite{CCR}
for details.
For a recent work dealing with the existence of weak solutions to a class of kinetic flocking models
of Cucker-Smale type we refer to \cite{KMT}.

\subsection*{Hydrodynamic description}
We provide the formal hydrodynamic limit for the system \eqref{kinetic}--\eqref{calG}.
Defining the mass and momentum densities 
\[
   \rho(t,x) = \int_{\R^d} f(t,x,v)\d v,\qquad
   \rho u(t,x) = \int_{\R^d} f(t,x,v) v\d v,
\]
we have
\[
   \td[\rho](x,y) &=& \int_{\R^d} \chi_{[0,1)}\left(\frac{|z-x|}{|y-x|}\right) \rho(z) \d z,\\
   \int_{\R^d} \calG[f](x,v) f(x,v) \d v &=& - \frac{1}{\gamma} \int_{\R^d} g\bigl(\td[\rho](x,y)\bigr) \rho(x)\rho(y) [u(x)-u(y)] \d y.
\]
The latter expression can be further simplified by noticing that the relative separation $\td[\rho](x,y)$
depends on $y$ only through $|x-y|$, so it can be written as $\td[\rho](x,y) = a_x(|x-y|)$, with
\[
    a_x(r) := \int_0^r R_x(s)\d s,\qquad R_x(s) := \int_{S_s(x)} \rho(y) \d S(y),
\]
where $\d S$ is the surface measure on the sphere $S_s(x)$ with radius $s>0$, centered at $x$.
Therefore, we have $\tot{}{r} a_x(r) = R_x(r)$ for every fixed $x\in\R^d$ and
\[
    \int_{\R^d} g\bigl(\td[\rho](x,y)\bigr) \rho(y) \d y = \int_0^\infty g\bigl(a_x(r)\bigr) \tot{}{r} a_x(r) \d r = 
     \int_0^1 g(s) \d s = \gamma,
\]
where we used the normalization $\int_{\R^d} \rho(x)\d x = a_x(\infty) = 1$.
Thus, we obtain
\(    \label{baru}
    \int_{\R^d} \calG[f](x,v) f(x,v) \d v = \rho(x) [\bar u(x) - u(x)]\qquad\mbox{with}\quad
    \bar u(x) = \frac{1}{\gamma} \int_{\R^d} g\bigl(\td[\rho](x,y)\bigr) \rho(y) u(y) \d y.
\)
Then, integrating the kinetic equation \eqref{kinetic}, \eqref{calG} againts the moments $(1,v)$
and taking the usual monokinetic closure assumption $f(t,x,v) = \rho(t,x)\delta(v-u(t,x))$,
we obtain the closed Euler-type system
\(
    \part{\rho}{t} + \grad\cdot(\rho u) &=& 0,  \label{Euler1} \\
    \part{u}{t} + (u\cdot\grad)u &=& \bar u - u. \label{Euler2}
\)
The right-hand side in the $u$-equation decribes the tendency of agents with velocity $u$
to relax to the local average velocity $\bar u$, defined in \eqref{baru}.
The word ``local'' is to be understood in the sense that
every agent measures the average velocity of a certain portion of its closest neighbors,
but independently of their actual metric distance.

The Euler system \eqref{Euler1}--\eqref{Euler2} falls into the class studied in Section 6 of \cite{Motsch-Tadmor}.
The authors adapt their method of active sets, originally developed for the discrete Cucker-Smale model,
to the continuum description. Their result on asymptotic flocking can be easily adapted to our situation under the assumption
that the communication rate function $g$ is strictly positive:

\begin{proposition}
Let $g(s) \geq g_0 > 0$ for all $s\in [0,1]$.
Consider the system \eqref{Euler1}--\eqref{Euler2} subject to compactly
supported inital data $(\rho^0,u^0)$ in $\R^d\times\R^d$ and assume that it admits global smooth
and compactly supported solutions $(\rho(t),u(t))$.
Then the position and, resp., velocity diameters of the solution $(\rho(t),u(t))$,
\[
   d_x(t) &=& \sup\{|x-y|,\; x,y\in\supp\rho(t)\},\\
   d_u(t) &=& \sup\{|u(t,x)-u(t,y)|,\; x,y\in\supp\rho(t)\},
\]
satisfy
\[
    \sup_{t\geq 0} d_x(t) < +\infty,\qquad \mbox{and}\quad d_u(t) \leq d_u(0) e^{-g_0^2 t}.
\]
\end{proposition}

\begproof
A slight modification of the proof of Proposition 6.4 in \cite{Motsch-Tadmor}.
\endproof

As already mentioned in Section \ref{sec:Formulation},
the assumption that $g(s) \geq g_0 > 0$ excludes
the interesting case when there is no \emph{direct} communication between certain
regions of $\supp\rho$, i.e., $g(\td[\rho](x,y))=0$ for $x\in\Omega_1$, $y\in\Omega_2$,
with some disjoint $\Omega_1$, $\Omega_2 \subset\supp\rho$.
This immediately leads to the question if it is possible to modify the method developed in Section \ref{sec:Flocking}
to prove asymptotic flocking for the hydrodynamic system \eqref{Euler1}--\eqref{Euler2}.
The natural way would be to approximate the solutions $(\rho,u)$
by finite systems of particles moving along characteristics, apply Theorem \ref{thm:flocking} for them,
and pass to the limit, using appropriate stability properties.
However, this program would be based on the a-priori assumption that the solution $(\rho,u)$
is uniformly approximable by strongly connected systems of particles (in the sense of Theorem \ref{thm:flocking}),
which is a very strong and technical assumption.
Therefore, we do not make any attempts in this direction.

\section{Extension to other models of collective dynamics}\label{sec:OtherModels}
The idea of introducing the discreet \eqref{topdist} and continuum \eqref{topdist-cont}
relative separation is of course applicable to the full spectrum of models of collective behavior.
Of course, this needs to be well justified from the modeling point of view,
as in many cases it may be more appropriate to stay with the metric interactions.

Let us provide just one particular example, the well-known attraction-repulsion particle model proposed
in \cite{attr-rep},
\(
    \dot x_i &=& v_i,  \label{AR1}\\
    \dot v_i &=& (a-b|v_i|^2)v_i - \frac{1}{N} \sum_{j\neq i} \grad U(|x_i-x_j|),\qquad i=1,\dots,N,  \label{AR2}
\)
where $a$ and $b$ are nonnegative parameters and $U:\R^d\to\R$ is a given potential modeling
the short-range repulsion and long-range attraction. The term corresponding
to $a$ models the self-propulsion of individuals, whereas the term corresponding to
$b$ is the friction assumed to follow Rayleigh's law. The balance of these two terms
imposes an asymptotic speed to the agent (if other effects are ignored), but does
not influence the orientation vector. A typical choice for $U$ is the Morse potential,
which is radial and given by
\(  \label{Morse}
   U(x) = - C_{A} e^{-|x|/\ell_A} + C_R e^{-|x|/\ell_R},
\)
where $C_A$, $C_R$ and $\ell_A$, $\ell_R$ are the strengths and the typical lengths of attraction and
repulsion, respectively.

Clearly, the pairwise forces of attraction and repulsion decay with the metric distance
between the respective agents. Let us now review the modeling assumptions
of these two effects in the context of biological collective behavior:
The repulsion force aims to describe the preference of the individuals to avoid
collisions or even close contact (for instance in case of canibalism) with conspecifics.
Therefore, it is a short-range metric interaction. On the other hand, the attraction
force aims to model the preference of the individuals to form a compact group,
e.g., a school, flock or swarm. The typical motivation for this behavior is
to avoid predation \cite{Vicsek-survey}. Therefore, the strength of the attractive
force should not be primarily related to the actual physical distance between individuals,
and it seems appropriate to introduce the relative separation to model the attraction.
For this we split the interaction term in \eqref{AR2} into two parts: the repulsion,
still modeled as a gradient of the repulsion potential
\(  \label{Ur}
   U_R(x) = C_R e^{-|x|/\ell_R},
\)
and the attraction, which, however, is not a potential gradient any more.
Instead, we propose the new form
\(   \label{AR2new}
    \dot v_i &=& (a-b|v_i|^2)v_i - \frac{1}{N} \sum_{j\neq i} \grad U_R(|x_i-x_j|) - \frac1N \sum_{j\neq i} g_A(\td(x_i,x_j)) \frac{x_i-x_j}{|x_i-x_j|},
    \qquad i=1,\dots,N,
\)
where $g_A:[0,1]\to[0,\infty)$ describes the strength of the attractive force in dependence
on the relative separation $\td(x_i,x_j)$ given by \eqref{topdist}.
The formal kinetic limit of \eqref{AR1}, \eqref{AR2new} reads then
\[
    \part{f}{t} + v\cdot\grad_x f + \grad_v\cdot[(a-b|v|^2)vf] - \grad_v\cdot[(\grad_x U_R\ast\rho)f]
    - \grad_v\cdot(\calG[f]f) = 0,
\]
with $\rho(t,x) = \int_{\R^d} f(t,x,v)\d v$ and
\[
   \calG[f](t,x) = \int_{\R^d} g(\td[f](x,y)) \frac{x-y}{|x-y|}\rho(t,y)\d y \,,
\]
with $\td[f](x,y)$ given by \eqref{topdist-cont}.

\subsection*{Numerical simulation}
\begin{figure}
{\centering \begin{tabular}[h]{cc}
\resizebox*{0.49\linewidth}{!}{\includegraphics{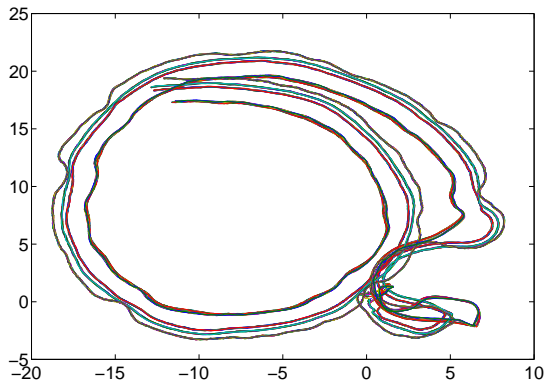}} &
\resizebox*{0.49\linewidth}{!}{\includegraphics{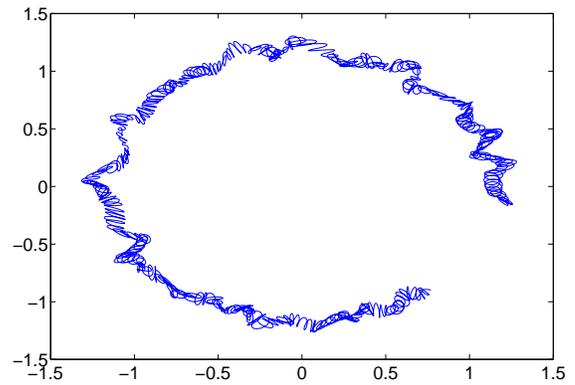}} \\
(a) all agents for $t\in[0,10]$ & (b) single agent for $t\in[8.0,8.5]$ \\
\resizebox*{0.49\linewidth}{!}{\includegraphics{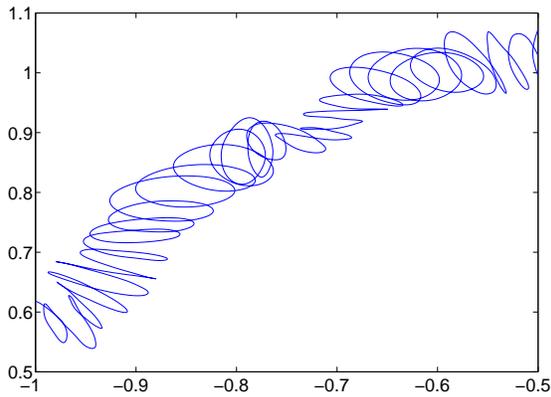}} &
\resizebox*{0.49\linewidth}{!}{\includegraphics{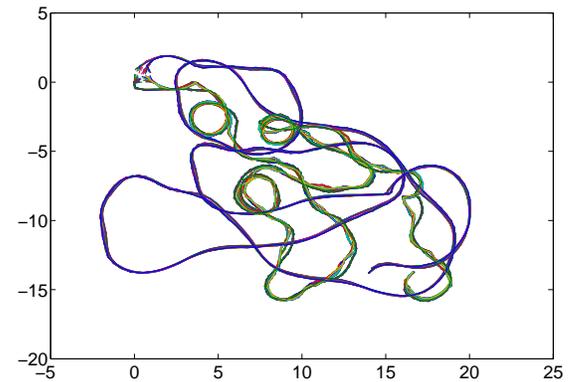}} \\
(c) single agents for $t\in[8.2,8.3]$ & (d) new simulation for $t\in[0,100]$ \\
\end{tabular}\par}
\caption{Simulation of the model \eqref{AR1}, \eqref{AR2new} with $N=100$ agents subject to a random initial condition
in positions and zero initial velocities.
(a) Upper left panel shows the trajectories of all agents during the time interval $t\in[0,10]$.
(b) A single agent trajectory for $t\in[8.0,8.5]$ is shown in the upper right panel,
with (c) another close-up in the lower left panel.
Finally, (d) another run of the simulation with a new random initial condition during the time interval $t\in[0,100]$ generated the pattern
shown in lower right panel, where again trajectories of all agents are recorded.}
\label{fig:newmodel}
\end{figure}

It is well known that the self-propelled attraction-repulsion model \eqref{AR1}--\eqref{AR2}
produces a variety of different patterns, called \emph{clumps, ring clumping,
rings} and \emph{mills} in \cite{attr-rep}. It is out of scope if this paper
to thoroughly investigate the patterning behavior of the  modified model \eqref{AR1}, \eqref{AR2new}.
Let us only mention that our preliminary numerical simulations suggest that
the new model is not only capable of reproducing all the above mentioned patterns,
but also exhibits interesting novel patterning phenomena.
In particular, we ran a simulation with $N=100$ agents,
initially randomly placed into the box $(0,1)\times(0,1)\subset\R^2$ with zero initial velocities.
We used the repulsion potential \eqref{Ur} with $C_R=1.0$ and $\ell_R=0.5$.
The attraction in \eqref{AR2new} was modeled using $g_A(s) = \frac{C_A}{\ell_A} e^{-s/\ell_a}$ with $C_A=1.0$
and $\ell_A=10/N=0.1$; this means that on average, each agent is mostly attracted by its 10 closest neighbors.
Note that this scaling ensures that the repulsive and attractive forces are initially, i.e.,
when all the agents are randomly distributed in the unit box, of the same order of magnitude.
We ran the simulation in the time interval $t\in[0,10]$ and recorded the trajectories of all the agents 
in the upper left panel of Fig. \ref{fig:newmodel}.
The agents quickly formed four groups of approximately same sizes and after an initial unordered movement,
the groups started to follow an approximately circular path, being imposed by the interaction between the groups.
However, on the scale of interactions within each group, we observe helixoidal-like movement of each single agent,
as depicted in the upper right panel of Fig. \ref{fig:newmodel} for the time interval $t\in[8.0,8.5]$
and another close-up in the lower left panel. It seems that this two-scale dynamic
is genuine to the model with topological attractive interactions and as far as we know, cannot be produced
by the original model \eqref{AR1}--\eqref{AR2}. Interestingly, the quasi-circular pattern on the scale of the groups
may not be persistent and may collapse into a chaotic movement, such as the one shown in the lower right panel of Fig. \ref{fig:newmodel}.
This last figure was produced by running the simulation with all the same parameters, starting from
another random sample for the initial condition, again with zero initial velocities.
We conclude that the modified model \eqref{AR1}, \eqref{AR2new} exhibits very interesting and new patterning dynamics,
which deserve further investigations. This is however not an objective of the present paper.

\noindent

\section{Acknowledgment}
The author acknowledges several interesting discussions and valuable
hints provided by Jan Vyb\'iral (TU Berlin),
Jos\'e Antonio Ca\~nizo (U. Birmingham) and Jos\'e Antonio Carrillo (Imperial College London).

\end{document}